\pdfoutput=1
\documentclass[12pt]{article}

\usepackage[font=footnotesize,margin=1cm]{caption}

\usepackage[margin=1in]{geometry}

\usepackage{setspace}
\usepackage[utf8]{inputenc} 
\usepackage[T1]{fontenc}    
\usepackage{hyperref}       
\usepackage{url}            
\usepackage{booktabs}       
\usepackage{amsfonts}       
\usepackage{nicefrac}       
\usepackage{microtype}      

\usepackage{amsmath,amssymb,bbm,stmaryrd,mathrsfs,url,amsthm}
\newtheorem{theorem}{Theorem}
\newtheorem{lemma}{Lemma}

\usepackage{amstext}
\usepackage{array,xcolor}
\usepackage{color,soul}
\usepackage{hyperref}
\usepackage{blkarray}
\usepackage{amsmath}

\usepackage{amsmath}
\usepackage[all,cmtip]{xy}
\usepackage{booktabs,makecell}

\usepackage[T1]{fontenc}
\usepackage[utf8]{inputenc}
\usepackage{multirow,bigdelim}
\usepackage{makeidx}
\usepackage{overpic}
\makeindex

\usepackage[export]{adjustbox}
\usepackage{graphicx,floatrow}
\usepackage{subcaption}

\DeclareMathOperator*{\argmin}{arg\,min}
\newcommand{\balpha}{{\boldsymbol{\alpha}}}
\newcommand{\bbeta}{{\boldsymbol{\beta}}}
\newcommand{\bxi}{{\boldsymbol{\xi}}}

\newcommand{\R}{\mathbb{R}}

\newcommand{\E}{\operatorname{E}}

\newcommand{\vct}[1]{\boldsymbol{#1}}
\newcommand{\mtx}[1]{\boldsymbol{#1}}

\newcommand{\<}{\langle}
\renewcommand{\>}{\rangle}

\newcommand{\set}[1]{\mathcal{#1}}

\DeclareMathOperator*{\minimize}{\text{minimize}}

\newcommand{\va}{\vct{a}}
\newcommand{\vb}{\vct{b}}
\newcommand{\vc}{\vct{c}}

\newcommand{\vg}{\vct{g}}
\newcommand{\vh}{\vct{h}}

\newcommand{\vm}{\vct{m}}

\newcommand{\vs}{\vct{s}}
\newcommand{\vt}{\vct{t}}
\newcommand{\vu}{\vct{u}}
\newcommand{\vv}{\vct{v}}
\newcommand{\vw}{\vct{w}}
\newcommand{\vx}{\vct{x}}
\newcommand{\vy}{\vct{y}}
\newcommand{\vz}{\vct{z}}

\newcommand{\vbeta}{\vct{\beta}}

\newcommand{\vxi}{\vct{\xi}}

\newcommand{\mA}{\mtx{A}}
\newcommand{\mB}{\mtx{B}}
\newcommand{\mC}{\mtx{C}}
\newcommand{\mD}{\mtx{D}}
\newcommand{\mE}{\mtx{E}}
\newcommand{\mF}{\mtx{F}}

\newcommand{\mI}{\mtx{I}}

\newcommand{\mP}{\mtx{P}}
\newcommand{\mQ}{\mtx{Q}}

\newcommand{\setA}{\set{A}}
\newcommand{\setB}{\set{B}}
\newcommand{\setC}{\set{C}}
\newcommand{\setD}{\set{D}}

\newcommand{\setN}{\set{N}}

\newcommand{\setR}{\set{R}}
\newcommand{\setS}{\set{S}}

\newcommand{\yl}{y_\ell}
\renewcommand{\sl}{s_\ell}
\newcommand{\wl}{w_\ell}
\newcommand{\xl}{x_\ell}

\newcommand{\al}{\va_\ell}
\newcommand{\bl}{\vb_\ell}

\newcommand{\cl}{\vc_\ell}

\newcommand{\mtil}{\tilde{\vm}}
\newcommand{\htil}{\tilde{\vh}}

\newcommand{\dm}{\delta \vm}
\renewcommand{\dh}{\delta \vh}

\DeclareMathOperator{\sign}{sign}

\def\l{\ell}

\newcommand{\ho}{\vh^\natural}

\newcommand{\mo}{\vm^\natural}

\newcommand{\wo}{\vw^\natural}
\newcommand{\xo}{\vx^{\natural}}

\newcommand{\alt}{\al^\intercal}
\newcommand{\blt}{\bl^\intercal}
\newcommand{\clt}{\cl^\intercal}

\newcommand{\tildeh}{\tilde{\vh}}
\newcommand{\tildem}{\tilde{\vm}}

\newcommand{\hath}{\hat{\vh}}
\newcommand{\hatm}{\hat{\vm}}

\newcommand{\PP}{\mathbb{P}}

\title{A convex program for bilinear inversion of sparse vectors}

\author{
  Alireza Aghasi\thanks{
  aaghasi@gsu.edu, J. Mack Robinson College of Business, GSU}, Ali Ahmed\thanks{ali.ahmed@itu.edu.pk, Department of Electrical Engineering, ITU, Lahore}, Paul Hand\thanks{p.hand@northeastern.edu, Department of Mathematics and College of Computer and Information Science, Northeastern University}\ \ and Babhru Joshi\thanks{babhru.joshi@rice.edu, Department of Computational and Applied Mathematics, Rice University}
}
\begin{document}

\maketitle
\begin{abstract}
  We consider the bilinear inverse problem of recovering two vectors, $\vx\in\R^L$ and $\vw\in\R^L$, from their entrywise product. We consider the case where $\vx$ and $\vw$ have known signs and are sparse with respect to known dictionaries of size $K$ and $N$, respectively.  Here,  $K$ and $N$ may be larger than, smaller than, or equal to $L$.  We introduce $\ell_1$-BranchHull, which is a convex program posed in the natural parameter space and does not require an approximate solution or initialization in order to be stated or solved. We study the case where $\vx$ and $\vw$ are $S_1$- and $S_2$-sparse with respect to a random dictionary and present a recovery guarantee that only depends on the number of measurements as {$L\geq\Omega(S_1+S_2)\log^{2}(K+N)$}. Numerical experiments verify that the scaling constant in the theorem is not too large.  One application of this problem is the sweep distortion removal task in dielectric imaging, where one of the signals is a nonnegative reflectivity, and the other signal lives in a known subspace, for example that given by dominant wavelet coefficients. We also introduce a variants of $\ell_1$-BranchHull for the purposes of tolerating noise and outliers, and for the purpose of recovering piecewise constant signals.  We provide an ADMM implementation of these variants and show they can extract piecewise constant behavior from real images.
\end{abstract}

\section{Introduction}

We study the problem of recovering two unknown signals $\vx$ and $\vw$ in $\R^{L}$ from observations $\vy =\mathcal{A}(\vw,\vx)$, where $\mathcal{A}$ is a bilinear operator. Let $\mB \in \R^{L\times K}$ and $\mC \in \R^{L\times N}$  such that $\vw = \mB\vh$ and $\vx =\mC\vm$ with $\|\vh\|_0 \leq S_1$ and $\|\vm\|_0\leq S_2$. Let the bilinear operator $\mathcal{A}:\R^L \times \R^L\rightarrow \R^L$ satisfy
\begin{equation}\label{eq:measurements}
	\vy = \mathcal{A}(\vw,\vx)=\vw \odot \vx,
\end{equation}

where $\odot$ denotes entrywise product. The bilinear inverse problem (BIP) we consider is to find $\vw$ and $\vx$ from $\vy$, $\mB$, $\mC$ and $\sign{(\vw)}$, up to the inherent scaling ambiguity. 

BIPs, in general, have many  applications in signal processing  and machine learning and include fundamental practical problems like phase retrieval \cite{fienup1982phase, candes2012solving, candes2013phaselift}, blind deconvolution \cite{ahmed2012blind, stockham1975blind, kundur1996blind,aghasi2016sweep}, non-negative matrix factorization \cite{hoyer2004non,lee2001algorithms}, self-calibration \cite{ling2015self}, blind source separation \cite{Gardy05source}, dictionary learning \cite{tosic2011dictionary}, etc. These problems are in general challenging and suffer from identifiability issues that make the solution set non-unique and non-convex. A common identifiability issue, also shared by the BIP in \eqref{eq:measurements}, is the scaling ambiguity. In particular, if $(\wo,\xo)$ solves a BIP, then so does $(c\wo, c^{-1}\xo)$ for any nonzero $c \in \R$. In this paper, we resolve this scaling ambiguity by finding the point in the solution set closest to the origin with respect to the $\ell_1$ norm.
\begin{figure}[H]
\floatsetup[subfigure]{captionskip=-20pt} 
\ffigbox[]{\hspace{-60pt}
    \begin{subfloatrow}[2]
      \ffigbox[]{%
      \begin{overpic}[width=0.4\textwidth,height=0.4\textwidth,tics=1]{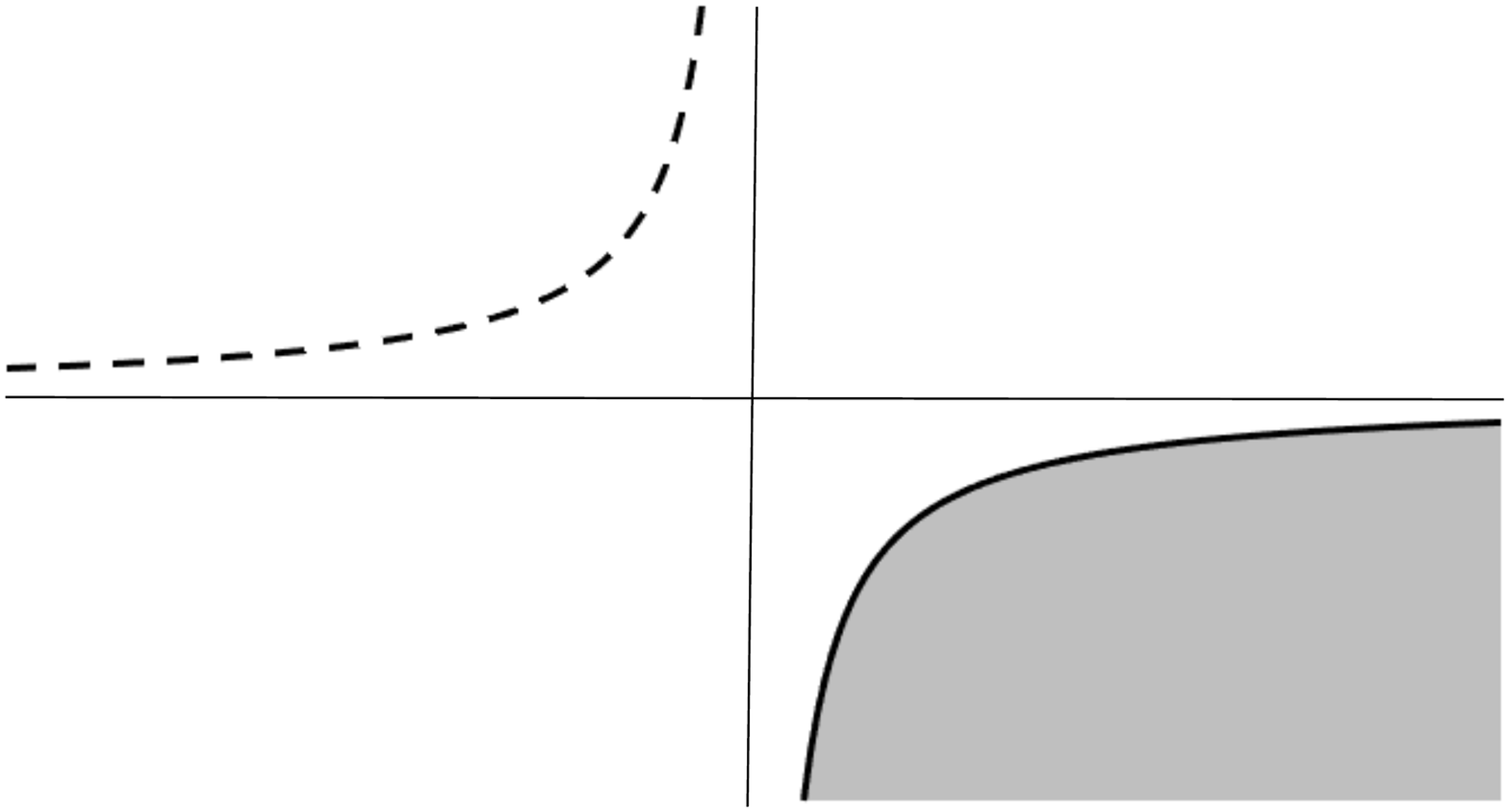}
     	\put(52,53){$0$}
		\put(49,80){$w_\ell$}
		\put(93,50){$x_\ell$}
		\put(45,34){\rotatebox{35}{$x_\ell w_\ell =y_\ell $}}
		\put(58,28.5){\rotatebox{0}{\scalebox{.93}{Convex Hull
		}}}
	  \end{overpic}
        }{\caption{Convex relaxation}\label{fig:hyperbola}}
      \hspace{\fill}\ffigbox[\FBwidth]{\raisebox{.8cm}{
        \begin{overpic}[scale = 0.3, trim= 8cm 11cm 6cm 4cm, clip = true]{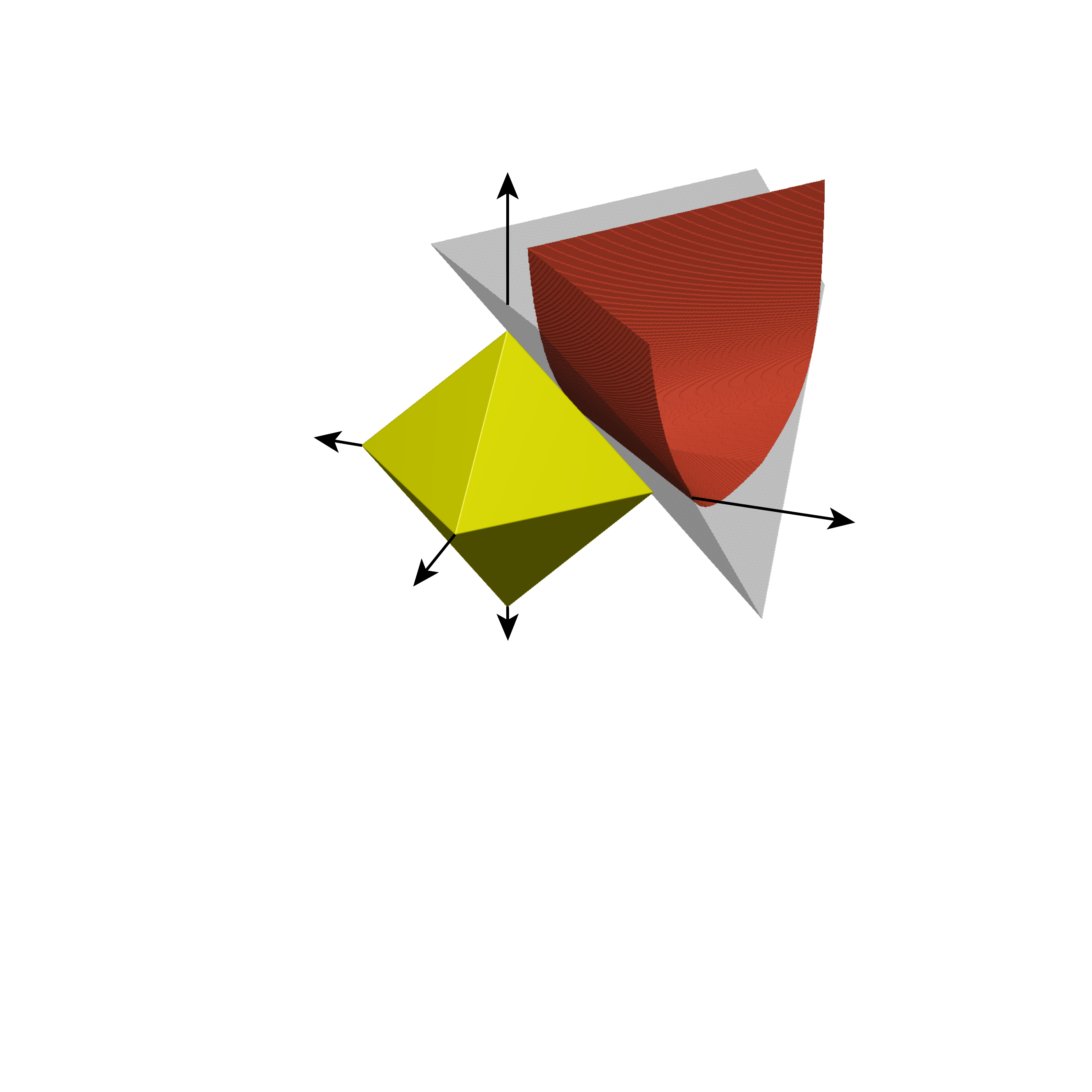}
		\put(16,8){$h_1$}
		\put(100,24){$h_2$}
		\put(33,92){$m_1$}
	\end{overpic}}}{\caption{Geometry of $\ell_1$-BranchHull}\label{fig:geometry}}
    \end{subfloatrow}
}{\caption{Panel (a) shows the convex hull of the relevant branch of a hyperbola given a measurement $\yl$ and the sign information $\sign(\wl)$. Panel (b) shows the interaction between the $\ell_1$-ball in the objective of $\eqref{eq:BH}$ with its feasibility set.  The feasibility set is `pointy' along a hyperbola, which allows for signal recovery where the $\ell_1$ ball touches it.   The gray hyperplane segments correspond to linearizations of the hyperbolic measurements, which is an important component of our recovery proof.   
}
	}
\end{figure}
\vspace{-10pt}
Another identifiability issue of the BIP in \eqref{eq:measurements} is if $(\wo,\xo)$ solves \eqref{eq:measurements}, then so does $(\bf{1}, \wo \odot \xo)$, where $\bf{1}$ is the vector of ones. {Unlike prior studies like \cite{ahmed2012blind}, where the signals are assumed to live in known subspaces, we resolve this structural ambiguity by additionally assuming the signals are sparse with respect to those known basis.  Natural choices for such bases include the standard basis, the Discrete Cosine Transform (DCT) basis, and a wavelet basis.}

Recent work on sparse BIP, specifically sparse blind deconvolution, in \cite{Bresler2016sparse} provides an exact recovery guarantee of the sparse vectors $\vh$ and $\vm$ that satisfy a "peakiness" condition, i.e. $\min\{\|\vh\|_\infty,\|\vm\|_\infty\}\geq c$ for some absolute constant $c \in \R$. This result holds with high probability for random measurements if the number of measurement, up to a log factor, satisfy $L \geq\Omega(S_1+S_2)$. For general vectors without the peakiness condition, the same work shows exact recover is possible if the number of measurements, up to a log factor, satisfy $L\geq \Omega(S_1S_2)$.

The main contribution of this paper is to introduce an algorithm for the sparse BIP described in \eqref{eq:measurements} that recovers the sparse vectors, which does not have to satisfy the peakiness condition {and does not need initialization}, under near optimal sample complexity. Precisely, we present a convex program stated in the natural parameter space, which in the noiseless setting with random $\mB$ and $\mC$, exactly recovers the sparse vectors with at most $S_1 + S_2$ combined nonzero entries with high probability if the number measurements satisfy $L \geq\Omega(S_1+S_2)\log^{2}(K+N)$.

\subsection{Convex program and main results}
We introduce a convex program written in the natural parameter space for the bilinear inverse problem described in \eqref{eq:measurements}. Let $(\ho,\mo) \in \R^{K}\times \R^{N}$ with $\|\ho\|_0 \leq S_1$ and $\|\mo\|_0\leq S_2$.  Let $w_\ell = \blt \ho$, $x_\ell = \clt \mo$ and $y_\ell = \blt\ho \clt\mo$, where $\blt$ and $\clt$ are the $\ell$th row of $\mB$ and $\mC$. Also, let $\vs = \sign(\vy)$ and $\vt = \sign(\mB\ho)$. The convex program we consider to recover $(\ho,\mo)$ is the $\ell_1$-BranchHull program
\begin{align}
\text{$\ell_1$-BH}: \qquad \underset{\vh\in \R^K, \vm \in \R^N}{\minimize}~\|\vh\|_1+\|\vm\|_1 \quad &\text{subject to}~~ \sl(\blt \vh \clt\vm) \geq |\yl|\label{eq:BH}\\
&\qquad\qquad\qquad t_\ell\blt\vh \geq 0, \quad \ell = 1,2,\ldots,L.\notag 
\end{align}
The motivation for the feasible set in program \eqref{eq:BH} follows from the observation that each measurement $y_\ell = w_\ell \cdot x_\ell$ defines a hyperbola in $\R^2$. As shown in Figure \eqref{fig:hyperbola}, the sign information $t_\ell = w_\ell$ restricts $(w_\ell,x_\ell)$ to one of the branch of the hyperbola. The feasible set in \eqref{eq:BH} corresponds to the convex hull of {particular branches of the hyperbola} for each $y_\ell$. This also implies that the feasible set is convex as it is the intersection of $L$ convex sets.

The objective function in \eqref{eq:BH} is an $\ell_1$ minimization over $(\vh,\vm)$ that finds a sparse point $(\hath, \hatm)$ with $\|\hath\|_1 = \|\hatm\|_1$. This scaling in the minimizer of \eqref{eq:BH} is justified by the observation that $(c\ho, c^{-1}\mo)$ is feasible for any non-zero $c \in \R$. So, the minimizer of \eqref{eq:BH}, under successful recovery, is $\left(\ho \sqrt{\frac{\|\mo\|_1}{\|\ho\|_1}},\mo \sqrt{\frac{\|\ho\|_1}{\|\mo\|_1}}\right)$. 

  Our main result is that under the structural assumptions that $\vw$ and $\vx$ live in random subspaces with $\ho$ and $\mo$ containing at most $S_1$ and $S_2$ non zero entries, the $\ell_1$-BranchHull program \eqref{eq:BH} recovers $\ho$, and $\mo$ (to within the scaling ambiguity) with high probability, provided  the number of measurements, up to log factors, satisfy $L \geq \Omega (S_1+S_2)\log^2(K+N)$.
\begin{theorem}\label{thm:Noiseless_Main}
	Suppose we observe the pointwise product of two vectors $\mB\ho$, and $\mC\mo$ through a bilinear measurement model in \eqref{eq:measurements}, where $\mB$, and $\mC$ are standard Gaussian random matrices. Then the $\ell_1$-BranchHull program \eqref{eq:BH} uniquely recovers $\left(\ho \sqrt{\frac{\|\mo\|_1}{\|\ho\|_1}},\mo \sqrt{\frac{\|\ho\|_1}{\|\mo\|_1}}\right)$ with probability at least $1-\mathrm{e}^{-(1/2)Lt^2}$ whenever  $L \geq 	C_t(S_1+S_2)\log^2(K+N)$, where $C_t$ is a constant that depends on $t \geq 0$. 
\end{theorem} 

\subsection{Prior art for bilinear inverse problems} 
Recent approaches to solving bilinear inverse problems like blind deconvolution and phase retrieval have been to lift the problems into a low rank matrix recovery task or to formulate an optimization programs in the natural parameter space. Lifting transforms the problem of recovering $\vh \in \R^K$ and $\vm \in \R^N$ from bilinear measurements to the problem of recovering a low rank matrix $\vh\vm^\intercal$ from linear measurements. The low rank matrix can then be recovered using a semidefinite program. The result in \cite{ahmed2012blind} for blind deconvolution showed that if $\vh$ and $\vm$ are representations of the target signals with respect to Fourier and Gaussian subspaces, respectively, then the lifting method successfully recovers the low rank matrix. The recovery occurs with high probability under near optimal sample complexity. Unfortunately, solving the semidefinite program is prohibitively computationally expensive because they operate in high-dimension space. Also, it is not clear how to enforce additional structure like sparsity of $\vh$ and $\vm$ in the lifted formulation in a way that allows optimal sample complexity \cite{li2013sparse, oymak2015simultaneously}. 

In comparison to the lifting approach for blind deconvolution and phase retrieval, methods that formulate an algorithm in the natural parameter space like alternating minimization and gradient descent based method are computationally efficient and also enjoy rigorous recovery guarantees under optimal or near optimal sample complexity \cite{li2016rapid, candes2014phase, netrapalli2013phase, sun2016geometric}. In fact, the work in \cite{Bresler2016sparse} for sparse blind deconvolution is based on alternating minimization. In the paper, the authors use an alternating minimization that successively approximate the sparse vectors while enforcing the low rank property of the lifted matrix. However, because these methods are non-convex, convergence to the global optimal requires a good initialization \cite{tu2015low, CC15, li2016rapid}. 

{Other approaches that operate in the natural parameter space include PhaseMax \cite{bahmani2016phase,  goldstein2016phasemax} and BranchHull \cite{aghasi2017branchHull}. PhaseMax is a linear program which has been proven to find the target signal in phase retrieval under optimal sample complexity if a good anchor vector is available. As with alternating minimization and gradient descent based approach, PhaseMax requires a good initialization. However, in PhaseMax the initialization is part of the optimization program  but in alternating minimization the initialization is part of the algorithmic implementation. BranchHull is a convex program which solves the BIP described in \eqref{eq:BH} excluding the sparsity assumption under optimal sample complexity. Like the $\ell_1$-BranchHull presented in this paper, BranchHull does not require an initialization but requires the sign information of the signals.}

 The $\ell_1$-BranchHull program \eqref{eq:BH} combines strengths of both the lifting method and the gradient descent based method. Specifically, the $\ell_1$-BranchHull program is a convex program that operates in the natural parameter space, {without a need for an initialization, and without restrictive assumptions on the class of recoverable signals.} These strengths are achieved at the cost of the sign information of the target signals $\vw$ and $\vx$. However, the sign assumption can be justified in imaging applications where the goal might be to recover pixel values of a target image, which are non-negative. Also, as in PhaseMax, the sign information can be thought of as an anchor vector which anchors the solution to one of the branches of the $L$ hyperbolic measurements.

\subsection{Extension to noise and outlier}
Extending the theory of the $\ell_1$-BranchHull program \eqref{eq:BH} to the case with noise is important as most real data contain significant noise.  {Formulation \ref{eq:BH} may be particularly susceptible to noise that changes the sign of even a single measurement.}  For the bilinear inverse problem as described in \eqref{eq:measurements} with small dense noise and arbitrary outliers, we propose the following robust $\ell_1$-BranchHull program
\begin{align}
	\text{RBH:}\quad \minimize_{\vh \in \R^K, \vm \in \R^N,\vxi \in \R^L} \|\vh\|_1+\|\vm\|_1+\lambda\|\vxi\|_{1}\quad \text{subject to }&s_\ell(\clt \vm + \xi_{\l})\blt \vh \geq |y_{\l}|,\label{eq:RBH}\\ 
	 &t_\ell \blt \vh \geq 0, \quad \l = 1,\dots, L \notag.
\end{align}	
The slack variable $\vxi$ controls the shape of the feasible set. For measurements $y_\ell$ with incorrect sign, the corresponding slack variables $\xi_\ell$ shifts the feasible set so that the target signal is feasible. In the outlier case, the $\l_1$ penalty promotes sparsity of slack variable $\vxi$. We implement a slight variation of the above program, detailed in Section \ref{total variation}, to remove distortions from real and synthetic images.

{
\subsection{Total variation extension of $\ell_1$-BranchHull}\label{total variation}
The robust $\ell_1$-BranchHull program \eqref{eq:RBH} is flexible and can be altered to remove distortions from an otherwise piecewise constant signal. In the case where $\vw = \mB\ho$ is a piecewise constant signal, $\vx = \mC\mo$ is a distortion signal and $\vy = \vw \odot\vx$ is the distorted signal, the total variation version \eqref{eq:TVBH} of the robust BranchHull program \eqref{eq:RBH}, under successful recovery, produces the piecewise constant signal $\mB\ho$, up to a scaling.
\begin{align}\label{eq:TVBH}
\text{TV BH}: \minimize_{\vh \in \R^K, \vm \in \R^N,\vxi \in \R^L}\hspace{-15pt}~\mbox{TV}\left(\mB\vh\right)+\|\vm\|_1+\lambda \|\boldsymbol{\xi}\|_1\quad \text{subject to}~~ &s_\ell(\xi_\ell+\vc_\ell^\top\vm)\vb_\ell ^\top \vh \geq |\yl|\\
&t_\ell\vb_\ell^\top \vh \geq 0, \quad \ell = 1,2,\ldots,L.\notag 
\end{align}
In \eqref{eq:TVBH}, $TV(\cdot)$ is a total variation operator and is the $\ell_1$ norm of the vector containing pairwise difference of neighboring elements of the target signal $\mB\vh$. We implement \eqref{eq:TVBH} to remove distortions from images in Section \ref{distortion removal}. 
}

\subsection{Notation}
Vectors and matrices are written with boldface, while scalars and entries of vectors are written in plain font.  For example, $c_\ell$ is the $\ell$the entry of the vector $\vc$.  We write $\boldsymbol{1}$ as the vector of all ones with dimensionality  appropriate for the context. We write $\mI_N$ as the $N\times N$ identity matrix. For any $x \in \R$, let $(x)_- \in \mathbb{Z}$ such that $x-1<(x)_-\leq x$. For any matrix $\mA$, let $\|\mA\|_F$ be the Frobenius norm of $\mA$. For any vector $\vx$, let $\|\vx\|_0$ be the number of non-zero entries in $\vx$. For $\vx \in \R^K$ and $\vy \in \R^N$, $(\vx,\vy)$ is the corresponding vector in $\R^K \times \R^N$, and $\<(\vx_1,\vy_1),(\vx_2,\vy_2)\> = \<\vx_1,\vx_2\> + \<\vy_1,\vy_2\>$. 

\section{Algorithm}\label{algorithm}
In this section, we present an Alternating Direction Method of Multipliers (ADMM) implementation of an extension of the robust $\ell_1$-BranchHull program \eqref{eq:RBH}. The ADMM implementation of the $\ell_1$-BranchHull program \eqref{eq:BH} is similar to the ADMM implementation of \eqref{eq:Rl1BH} and we leave it to the readers. The extension of the robust $\ell_1$-BranchHull program we consider is 
\begin{align}\label{eq:Rl1BH}
\underset{\vh \in \R^K, \vm \in \R^N,\vxi \in \R^L}{\minimize}~\|\mP\vh\|_1+\|\vm\|_1+\lambda \|\boldsymbol{\xi}\|_1\quad &\text{subject to}~~ s_\ell(\xi_\ell+\vc_\ell^\top\vm)\vb_\ell ^\top \vh \geq |\yl|\\
&\qquad\qquad t_\ell\vb_\ell^\top \vh \geq 0, \quad \ell = 1,2,\ldots,L,\notag 
\end{align}
where $\mP \in \R^{J \times K}$ for some $J \in \mathbb{Z}$. The above extension reduces to the robust $\ell_1$-BranchHull program if $\mP = \mI_K$. Recalling that $\vw = \mB\vh$ and $\vx = \mC\vm$, we make use of the following notations
\[\vu = \begin{pmatrix} \vx\\ \vw\\ \boldsymbol{\xi}\end{pmatrix}, ~~ \vv = \begin{pmatrix} \vm\\ \vh\\ \lambda\boldsymbol{\xi}\end{pmatrix} , ~~ \mE = \begin{pmatrix} \mC&\boldsymbol{0}&\boldsymbol{0}\\ \boldsymbol{0}& \mB&\boldsymbol{0} \\ \boldsymbol{0}& \boldsymbol{0}& \lambda^{-1}\mI_L\end{pmatrix}\mbox{and}~~ \mQ = \begin{pmatrix} \mI_N &\boldsymbol{0}&\boldsymbol{0}\\ \boldsymbol{0}& \mP &\boldsymbol{0} \\ \boldsymbol{0}& \boldsymbol{0}&\mI_L\end{pmatrix}.
\]
Using this notation, our convex program can be compactly written as
\begin{align}\underset{\vv\in \R^{N+K+L},\vu\in \R^{3L}}{\minimize}~~\|\mQ\vv\|_1~~\text{subject to}~~ \vu = \mE\vv ,~~ \vu\in \mathcal{C}. \notag
\end{align}
Here $C = \left\{(\vx,\vw,\vxi)\in\R^{3L}|\ s_\ell(\xi_\ell+x_\ell)w_\ell\geq |y_\ell|,\ t_\ell w_\ell\geq 0,\ \l = 1,\dots, L\right\}$ is the convex feasible set of $\eqref{eq:Rl1BH}$. Introducing a new variable $\vz$ the resulting convex program can be written as
\[\underset{\vv,\vu,\vz}{\minimize}~~\|\vv\|_1~~\text{subject to}~~ \vu = \mE\vz, ~~ \mQ\vz = \vv ,~~ \vu\in \mathcal{C}. \notag
\]
We may now form the scaled ADMM steps as follows
\begin{align}\label{admm1}
\vu_{k+1} &= \argmin_{\vu}~~ \mathcal{I}_{\mathcal{C}}(\vu) + \frac{\rho}{2}\left\|\vu + \balpha_k   - \mE\vz_k \right\|^2\\\label{admm2}
\vv_{k+1} & = \argmin_{\vv}~~ \|\vv\|_1 + \frac{\rho}{2}\left\|\vv + \bbeta_k  - \mQ\vz_k \right\|^2\\
\vz_{k+1} &= \argmin_{\vz}~~ \frac{\rho}{2}\left\|\balpha_k + \vu_{k+1} - \mE\mQ\vz \right\|^2 +  \frac{\rho}{2}\left\|\bbeta_k + \vv_{k+1} - \mQ\vz \right\|^2\label{admm3}\\ \notag \balpha_{k+1} &= \balpha_k + \vu_{k+1} - \mE\vz_{k+1},\\ \bbeta_{k+1} &= \bbeta_k + \vv_{k+1} - \mQ\vz_{k+1}.\notag
\end{align}
{where $\mathcal{I}_C(\cdot)$ in \eqref{admm1} is the indicator function on $C$ such that $\mathcal{I}_C(u) =0$ if $u \in C$ and infinity otherwise}. We would like to note that the first three steps of the proposed ADMM scheme can be presented in closed form. The update in \eqref{admm1} is the following projection
\[\vu_{k+1} = \mbox{proj}_{\mathcal{C}}\left(  \mE\vz_k - \balpha_k \right),
\]
where $\mbox{proj}_{\mathcal{C}}(\vz)$ is the projection of $\vz$ onto $\mathcal{C}$. Details of computing the projection onto C are presented in the Supplementary material. The update in \eqref{admm2} can be written in terms of the soft-thresholding operator
\[\vv_{k+1} = S_{1/\rho} \left( \mQ\vz_k - \bbeta_k \right) , \qquad\mbox{where}\quad \left(S_{c}(\vv)\right)_i = \left\{ \begin{array}{cc} v_i-c & v_i>c\\ 0 & |v_i|\leq c\\ v_i+c & v_i<-c\end{array}\right., 
\] 
 where $c>0$ and $\left(S_{c}(\vv)\right)_i$ is the $i$th entry of $S_{c}(\vv)$. Finally, the update in \eqref{admm3} takes the following form
\[\vz_{k+1} = \left(\mE^\top\mE + \mQ^\intercal\mQ \right)^{-1} \left(\mE^\top\left( \balpha_k + \vu_{k+1}\right) + \mQ^\intercal(\bbeta_k + \vv_{k+1}) \right).
\]
In our implementation of the ADMM scheme, we initialize the algorithm with the $\vz_0 = \bf{0}$, $\alpha_0 = \bf{0}$, $\vbeta_0 = \bf{0}$. 

\section{Numerical Experiments}
In this section, we provide numerical experiments on synthetic and real data. The synthetic experiment numerically verifies Theorem \ref{thm:Noiseless_Main} with a low scaling constant. The experiment on real data  shows total variation $\ell_1$-BranchHull program can be used to remove distortions from an image.

\subsection{Phase Portrait}
\begin{figure}[H]
\centering
\includegraphics[scale = 0.25]{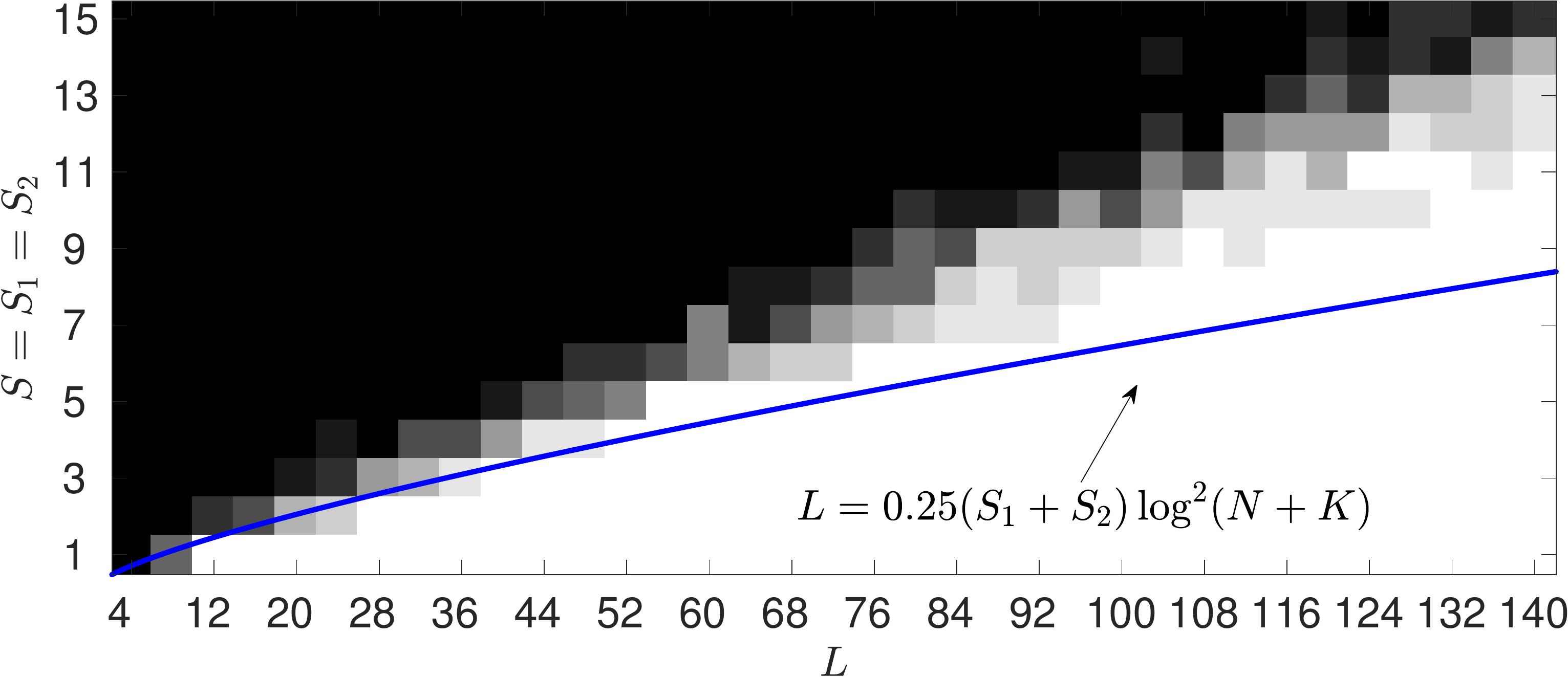}
\caption{The empirical recovery probability from synthetic data with sparsity  level $S$ as a function of total number of measurements $L$. Each block correspond to the average from 10 independent trials. White blocks correspond to successful recovery and black blocks correspond to unsuccessful recovery. The area to the right of the line satisfies $L > 0.25(S_1+S_2)\log^2(N+K)$.}
\label{phase_plot_noiseless}
\end{figure}
We first show a phase portrait that verifies Theorem \ref{thm:Noiseless_Main}. Consider the following measurements: fix $N \in \{20,40,\dots, 300\}$, $L \in \{4,8,\dots,140\}$ and let $K = N$. Let the target signal $(\ho, \mo) \in \mathbb{R}^{K}\times \mathbb{R}^{N}$ be such that both $\ho$ and $\mo$ have $0.05N$ non-zero entries with the nonzero indices randomly selected and set to $\pm 1$. Let $S_1$ and $S_2$ be the number of nonzero entries in $\ho$ and $\mo$, respectively. Let $\mB \in \mathbb{R}^{L\times K}$ and $\mC \in \mathbb{R}^{L\times N}$ such that $B_{ij}\sim \frac{1}{\sqrt{L}}\mathcal{N}(0,1)$ and $C_{ij}\sim \frac{1}{\sqrt{L}}\mathcal{N}(0,1)$. Lastly, let $\vy = \mB\ho\odot \mC\mo$ and $\vt = \text{sign}(\mB\ho)$.

Figure \ref{phase_plot_noiseless} shows the fraction of successful recoveries from 10 independent trials using \eqref{eq:BH} for the bilinear inverse problem \eqref{eq:measurements} from data as described above. Let $(\hat{\vh},\hat{\vm})$ be the output of \eqref{eq:BH} and let $(\tilde{\vh},\tilde{\vm})$ be the candidate minimizer. We solve \eqref{eq:BH} using an ADMM implementation similar to the ADMM implementation detailed in Section \ref{algorithm} with the step size parameter $\rho = 1$. For each trial, we say \eqref{eq:BH} successfully recovers the target signal if  $\|(\hat{\vh},\hat{\vm})-(\tilde{\vh},\tilde{\vm})\|_{2} <10^{-10}$. Black squares correspond to no successful recovery and white squares correspond to 100\% successful recovery. The line corresponds to $L = C(S_1+S_2)\log^{2}(K+N)$ with $C = 0.25$ and indicates that the sample complexity constant in Theorem \ref{thm:Noiseless_Main} is not very large. 

\subsection{Distortion removal from images}\label{distortion removal}
\begin{figure}[H]
\centering
	\begin{subfigure}{.24\textwidth}\centering
		\captionsetup{skip=5pt,oneside,margin={0cm,0cm}} 
		\includegraphics[scale = .16]{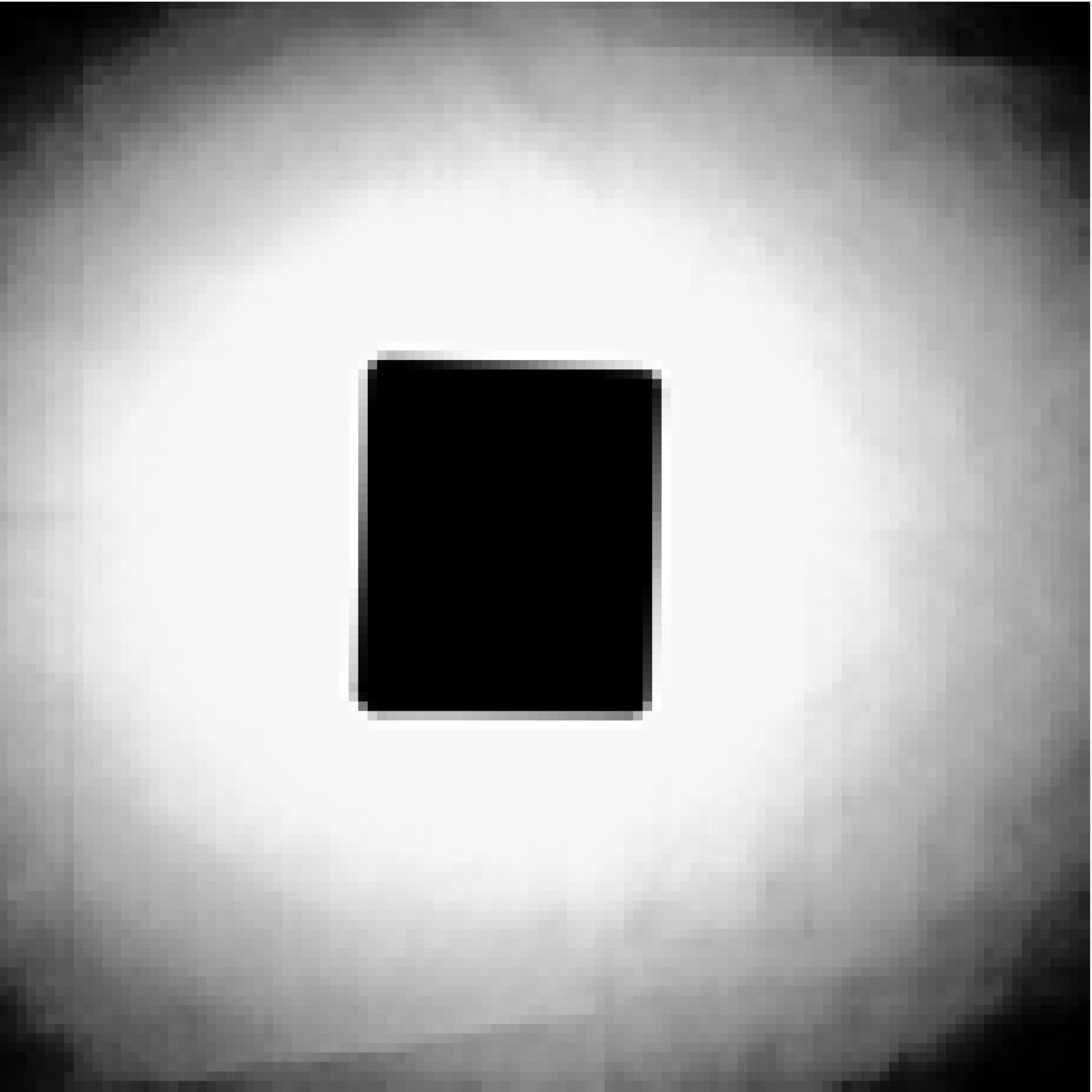}
		\caption{Distorted image}
		\label{data_y}
	\end{subfigure}
	\begin{subfigure}{.24\textwidth}\centering
		\captionsetup{skip=5pt,oneside,margin={-0cm,0cm}} 
		\centering
		\includegraphics[scale = .16,frame]{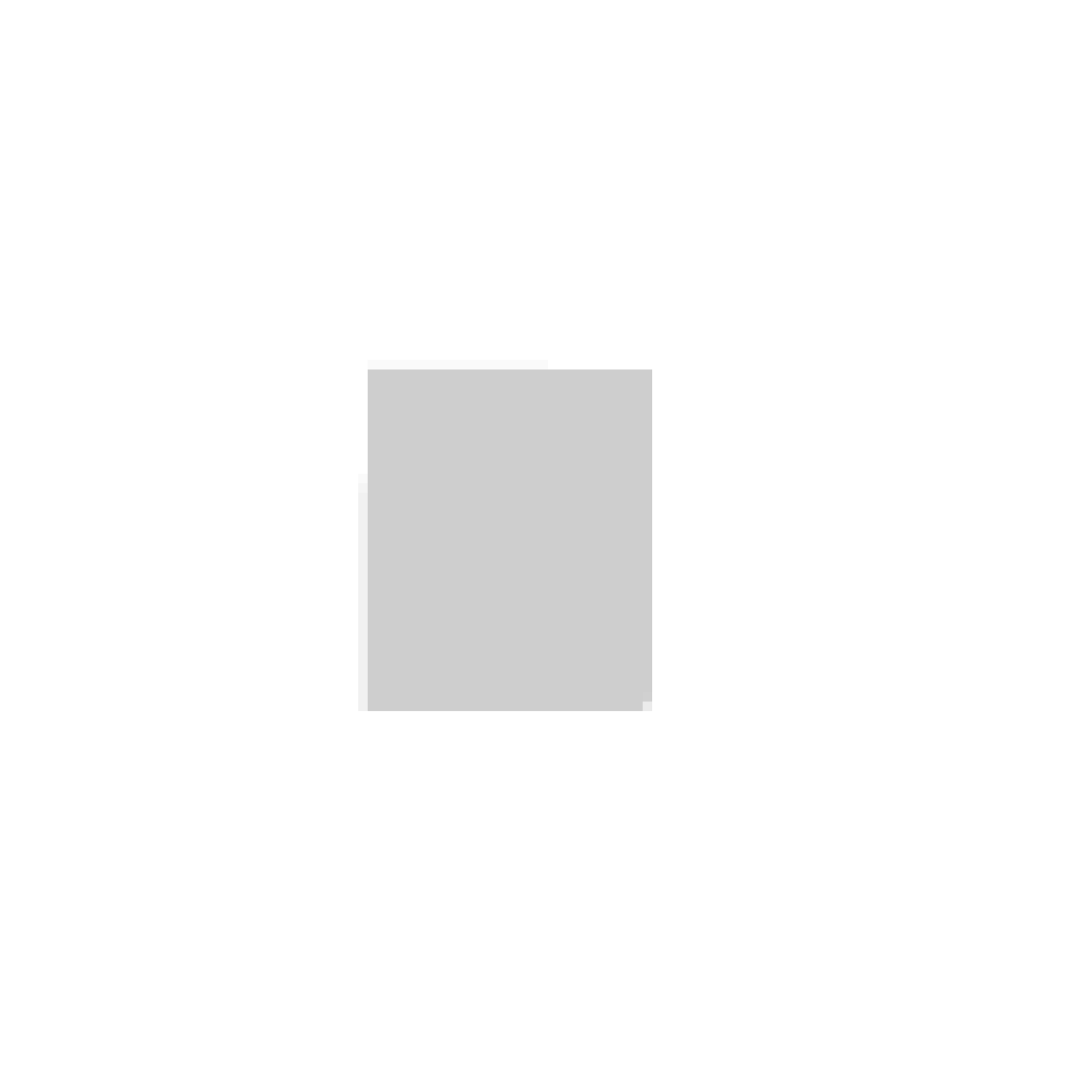} 
		\caption{Recovered image}
		\label{output_h}
	\end{subfigure}
	\begin{subfigure}{.24\textwidth}\centering
		\captionsetup{skip=5pt,oneside,margin={0cm,0cm}} 
		\includegraphics[scale = .16]{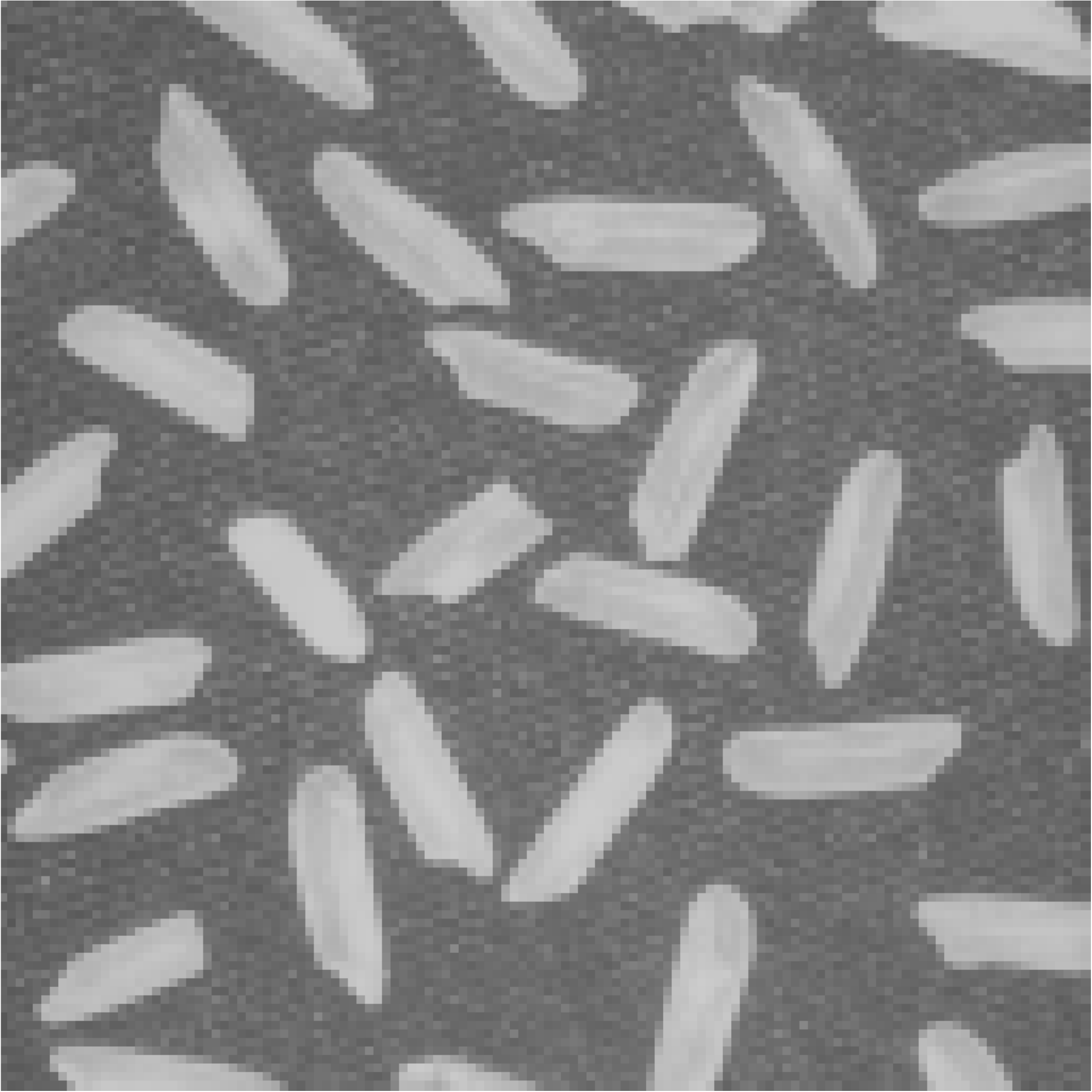}
		\caption{Distorted image}
		\label{rice_y}
	\end{subfigure}
	\begin{subfigure}{.24 \textwidth}\centering
		\captionsetup{skip=5pt,oneside,margin={-0cm,0cm}} 
		\centering
		\includegraphics[scale = .16]{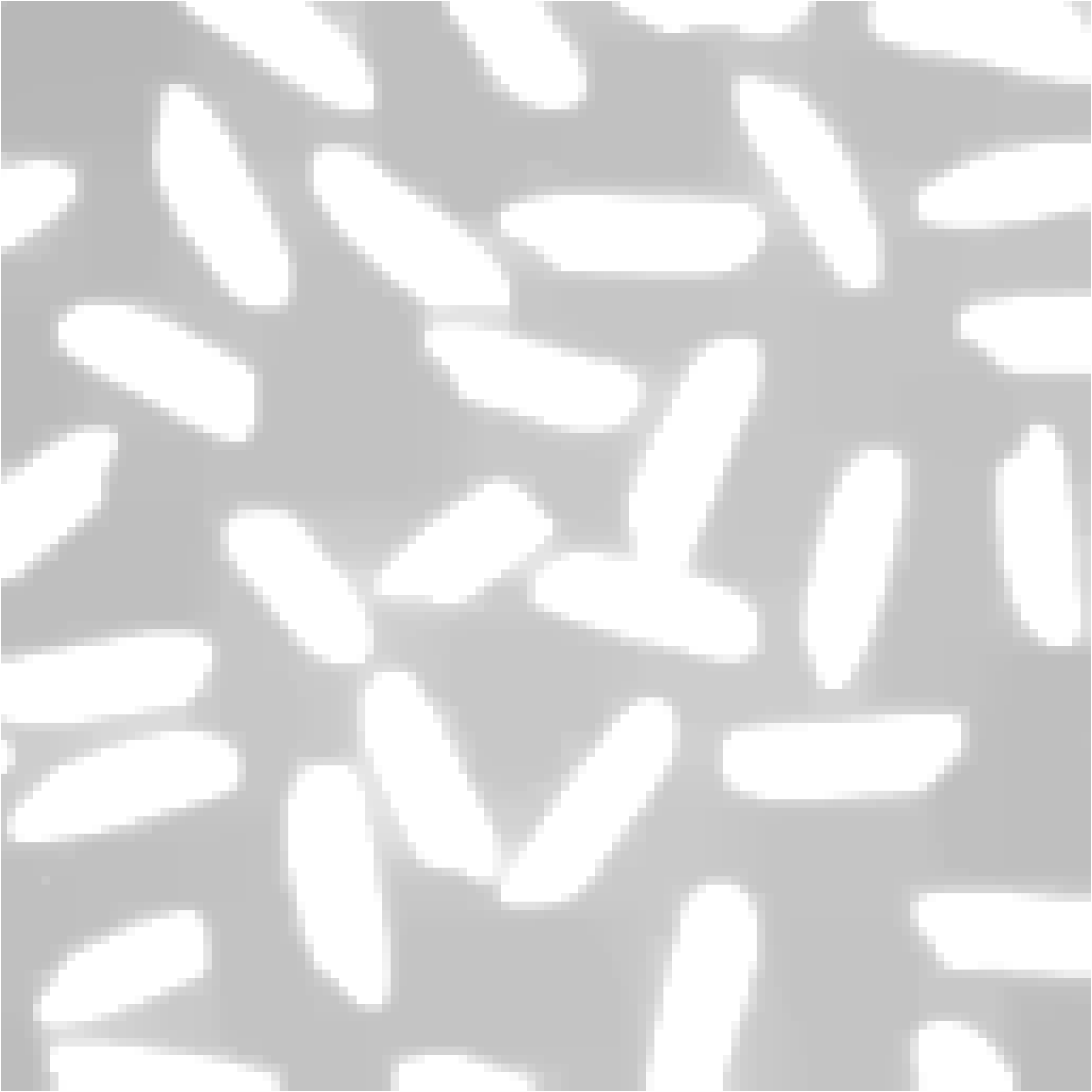} 
		\caption{Recovered image}
		\label{rice_h}
	\end{subfigure}\vspace{5pt}
	\caption{Panel (a) shows an image of a mousepad with distortions and panel(b) is the piecewise constant image recovered using total variation $\ell_1$-BranchHull. Similarly, panel (d) shows an image containing rice grains and panel (e) is the recovered image.}
\end{figure}
\vspace{-10pt}
{We use the total variation BranchHull program \eqref{eq:TVBH} to remove distortions from real images $\tilde{\vy} \in \R^{p\times q}$. In the experiments, The observation $\vy \in \mathbb{R}^L$ is the column-wise vectorization of the image $\tilde{\vy}$, the target signal $\vw=\mB\vh$ is the vectorization of the piecewise constant image and $\vx=\mC\vm$ corresponds to the distortions in the image. We use \eqref{eq:TVBH} to recover piecewise constant target images like in the foreground of Figure \ref{data_y} with $TV(\mB\vh) =\|\mD \mB\vh\|_1$, where $\mD = \left[\begin{array}{c}\mtx{D}_v \\ \mtx{D}_h\end{array}\right]$ in block form. Here, $\mD_v \in \mathbb{R}^{(L-q)\times L}$ and  $\mD_h \in \mathbb{R}^{(L-p)\times L}$ with
\[(\mD_v)_{ij} = \left\{\begin{array}{c l} -1 & \text{if } j=i+\left(\frac{i-1}{p-1}\right)_- \\ 1 & \text{if } j=i+1+\left(\frac{i-1}{p-1}\right)_-\\ 0 & \text{otherwise}\end{array} \right., \ (\mD_h)_{ij} = \left\{\begin{array}{c l} -1 & \text{if } j=i\\ 1 & \text{if } j=i+p\\ 0 & \text{otherwise}\end{array} \right. .
\]
}
Lastly, we solve \eqref{eq:TVBH} using the ADMM algorithm detailed in Section \ref{algorithm} with $\mP=\mD\mB$.

%

 We now show two experiments on real images. The first image, shown in Figure \ref{data_y}, was captured using a camera and resized to a $115 \times 115$ image. The measurement $\vy \in \mathbb{R}^{L}$ is the vectorization of the image with $L = 13225$. Let $\mB$ be the $L \times L$ identity matrix. Let $\mF$ be the $L \times L$ inverse DCT matrix. Let $C \in \mathbb{R}^{L \times 300}$ with the first column set to $\bf{1}$ and remaining columns randomly selected from columns of $\mF$ without replacement. The matrix $\mC$ is scaled so that $\|\mC\|_F = \|\mB\|_F = \sqrt{L}$. The vector of known sign $\vt$ is set to $\bf{1}$. Let $(\hat{\vh},\hat{\vm},\hat{\vxi})$ be the output of \eqref{eq:TVBH} with $\lambda = 10^3$ and $\rho = 10^{-4}$. Figure \ref{output_h} corresponds to $\mB\hat{\vh}$ and shows that the object in the center was successfully recovered.

The second real image, shown in Figure \ref{rice_y}, is an image of rice grains. The size of the image is $128 \times 128$. The measurement $\vy \in \mathbb{R}^{L}$ is the vectorization of the image with $L = 16384$. Let $B$ be the $L \times L$ identity matrix. Let $C \in \mathbb{R}^{L \times 50}$ with the first column set to $\bf{1}$. The remaining columns of $\mC$ are sampled from Bessel function of the first kind $J_{\nu}(z)$ with each column corresponding to a fixed $z$. Specifically, let $\vg \in \mathbb{R}^{L}$ with $g_i = -9+14\frac{i-1}{L-1}$. For each remaining column $\vc$ of $\mC$, let $\vz \sim \mathcal{N}(\mtx{0},I_3)$ and $c_i = J_{\frac{g_i}{6+0.1|z_1|}+5|z_2|}(0.1+10|z_3|)$. The matrix $\mC$ is scaled so that $\|\mC\|_F = \|\mB\|_F =\sqrt{L}$. The vector of known sign $\vt$ is set to $\bf{1}$. Let $(\hat{\vh},\hat{\vm},\hat{\vxi})$ be the output of \eqref{eq:TVBH} with $\lambda = 10^3$ and $\rho = 10^{-7}$. Figure \ref{rice_h} corresponds $\mB\hat{\vh}$. 

\section{Proof Outline}
In this section, we provide a proof of Theorem \ref{thm:Noiseless_Main} by considering a related linear program with larger feasible set.   Let $(\ho,\mo) \in   \R^{K}\times \R^{N}$ with $\|\ho\|_0 \leq S_1$ and $\|\mo\|_0\leq S_2$.  Let $w_\ell = \blt \ho$, $x_\ell = \clt \mo$ and $y_\ell = \blt\ho \cdot \clt\mo$. Also, let $\vs = \sign(\vy)$ and $\vt = \sign(\mB\ho)$. We will shows that the \eqref{eq:BH} recovers $(\tildeh,\tildem)$ such that $(\tildeh,\tildem)=\left(\ho \sqrt{\frac{\|\mo\|_1}{\|\ho\|_1}},\mo \sqrt{\frac{\|\ho\|_1}{\|\mo\|_1}}\right)$.

 Consider program \eqref{eq:LP} which has a linear constraint set that contains the feasible set of the $\ell_1$-BrachHull program \eqref{eq:BH}.
\begin{align}\label{eq:LP}
\text{LP}: \qquad \underset{\vh\in \R^K, \vm \in \R^N}{\minimize}~\|\vh\|_1+\|\vm\|_1\text{subject to}~~&\sl(\blt\vh\clt\tilde{\vm}+\blt\tilde{\vh}\clt\vm) \geq 2 |\yl|\\
&\ell = 1, 2, \ldots, L,\notag
\end{align}
Let 
\begin{align}\label{eq:setS}
\setS:= \left\{ (\vh,\vm) \in \R^K \times \R^N \ | \ (\vh,\vm)=  \alpha (-\tilde{\vh},\tilde{\vm}), \ \text{and} \ \alpha \in [-1,1] \right\}.
\end{align}
 Observe that if $(\tilde{\vh},\tilde{\vm})$ is a minimizer of \eqref{eq:LP} then so are  all the points $(\tilde{\vh},\tilde{\vm}) \oplus \setS$.
\begin{lemma}\label{lem:feas}
If the optimization program \eqref{eq:LP} recovers $(\vh,\vm) \in (\htil,\mtil)\oplus \setS $, then the BranchHull program \eqref{eq:BH} recovers $(\htil,\mtil)$.	
\end{lemma}
 A proof of Lemma \ref{lem:feas}{, provided in Supplementary material,} follows from the observations that the feasible set of \eqref{eq:LP} contains the feasible set of \eqref{eq:BH} and $(\htil,\mtil)$ is the only feasible point in \eqref{eq:BH} among all  $(\vh,\vm) \in (\htil,\mtil)\oplus \setS $.
 
  We now show that the only $(\vh,\vm) \in  (\tilde{\vh},\tilde{\vm}) \oplus \setS$ minimizes \eqref{eq:LP}. Let $\alt = (\clt\tilde{\vm}\blt,\blt\tilde{\vh}\clt) \in \R^{K+N}$ denote the $\ell$th row of a matrix $\mA$. The linear constraint in \eqref{eq:LP} are now simply $\vs \odot \mA(\vh,\vm)\geq 2|\vy|$. Note that $
\setS \subset \setN : = \text{span}(-\tilde{\vh},\tilde{\vm}) \subseteq \text{Null}(\mA)$.

Our strategy will be to show that for any feasible perturbation $(\dh,\dm) \in \setN_{\perp}$ the objective of the linear program \eqref{eq:LP} strictly increases, where $\setN_{\perp}$ is the orthogonal complement of the subspace $\setN$.  This will be equivalent to showing that the solution of \eqref{eq:LP} lies in the set $(\tildeh,\tildem)\oplus \setS$. 

The subgradient of the $\ell_1$-norm at the proposed solution $(\tilde{\vh},\tilde{\vm})$ is
\[
\partial \|(\tilde{\vh},\tilde{\vm})\|_1 := \{\vg\in \R^{K+N}: \|\vg\|_\infty \leq 1 ~\text{and}~ \vg_{\Gamma_h} = \sign(\ho_{\Gamma_h}) \ ,\vg_{\Gamma_m} = \text{sign}(\mo_{\Gamma_m})\},
\]
where $\Gamma_h,$ and $\Gamma_m$ denote the support of non-zeros in $\ho$, and $\mo$, respectively. To show the linear program converges to a solution $(\hat{\vh},\hat{\vm}) \in (\htil,\mtil)\oplus \mathcal{S}$, it suffices to show that the set of following descent directions
\begin{align}\label{eq:setD}
&\left\{ (\dh,\dm) \in \setN_{\perp}:\big\<\vg,(\dh,\dm)\big\> \leq 0,~ \forall \vg \in \partial \|(\tilde{\vh},\tilde{\vm})\|_1  \right\}\notag\\
&\subseteq  \left\{(\dh,\dm) \in \setN_{\perp}: \<\vg_{\Gamma_h},\dh_{\Gamma_h}\>+\<\vg_{\Gamma_m},\dm_{\Gamma_m}\>+\|(\dh_{\Gamma_h^c},\dm_{\Gamma_m^c})\|_1 \leq 0\right\}\notag\\
&\subseteq \left\{ (\dh,\dm) \in \setN_{\perp}: -\|\vg_{\Gamma_h \cup \Gamma_m}\|_2\|(\dh_{\Gamma_h},\dm_{\Gamma_m})\|_2+\|(\dh_{\Gamma_h^c},\dm_{\Gamma_m^c})\|_1\leq 0\right\}\notag\\
&= \left\{ (\dh,\dm) \in \setN_{\perp}: \|(\dh_{\Gamma_h^c},\dm_{\Gamma_m^c})\|_1 \leq \sqrt{S_1+S_2}\|(\dh_{\Gamma_h},\dm_{\Gamma_m})\|_2\right\} =: \mathcal{D}
\end{align}
does not contain any vector $(\dh,\dm)$ that is consistent with the constraints. We do this by quantifying the ``width" of the set $\setD$ through a Rademacher complexity, and a probability  the gradients of the constraint functions lie in a certain half space. This allows us to use small ball method \cite{koltchinskii2015bounding,mendelson2014learning} to ultimately show that it is highly unlikely to have descent directions in $\setD$ that meet the constraints in \eqref{eq:LP}. We now concretely state the definitions of the Rademacher complexity, and probability term mentioned above. 

Define linear functions 
\[
f_\ell(\vh,\vm) : = \left\<(\blt\tildeh\cl,\clt\tildem \bl), (\vh,\vm)\right\>, \ell = 1,2,3,\ldots, L.
\]
The linear constraints in the LP \eqref{eq:LP} can be expressed as $\sl f_\ell(\vh,\vm) {\geq} 2|y_\ell|$. The gradients of $f_\ell$ w.r.t. $(\vh,\vm)$ are then simply $\nabla f_\ell = (\tfrac{\partial f_\ell(\vh,\vm)}{\partial \vh}, \tfrac{\partial f_\ell(\vh,\vm)}{\partial \vm}) = (s_\ell\clt\tildem\bl,s_\ell\blt\tildeh\cl)$. Define the Rademacher complexity of a set $\setD \subset \R^{M}$ as
\begin{align}\label{eq:Rademacher-Complexity}
\mathfrak{C}(\setD) := \E\sup_{(\vh,\vm)\in \setD} \tfrac{1}{\sqrt{L}}\sum_{\ell=1}^L\varepsilon_\ell \left\<\nabla f_\ell,\tfrac{ (\vh,\vm)}{\|(\vh,\vm)\|_2}\right\>,
\end{align}
where $\varepsilon_1,\varepsilon_2, \ldots, \varepsilon_L$ are iid Rademacher random variables independent of everything else. For a set $\setD$, the quantity $\mathfrak{C}(\setD)$ is a measure of width of $\setD$ around the origin interms of the gradients of the constraint functions. For example, an 
equally distributed random set of gradient functions might lead to a smaller value of $\setC(\setD)$. 

Our results also depend on a probability $p_{\tau}(\setD)$, and a positive parameter $\tau$ introduced below
\begin{align}\label{eq:Tail-Prob}
\mathfrak{p}_{\tau}(\setD) = \inf_{(\vh,\vm)\in \mathcal{D}} \PP\left(\left\<\nabla f_\ell,\tfrac{ (\vh,\vm)}{\|(\vh,\vm)\|_2}\right\>\geq \tau\right).
\end{align}
Intuitively, $p_{\tau}(\setD)$ quantifies the size of $\setD$ through the gradient vectors. For a small enough fixed parameter, a small value of $p_{\tau}(\setD)$ means that the $\setD$ is mainly invisible to to the gradient vectors. 

\begin{lemma}\label{lem:Mendelson}
	Let $\setD$ be the set of descent directions, already characterized in \eqref{eq:setD}, for which $\mathfrak{C}(\setD)$, and $\mathfrak{p}_{\tau}(\setD)$  can be determined using  \eqref{eq:Rademacher-Complexity}, and \eqref{eq:Tail-Prob}. Choose $ L \geq \left(\frac{2\mathfrak{C}(\mathcal{D})+t\tau}{\tau \mathfrak{p}_{\tau}(\mathcal{D})}\right)^2$ for any $t >0$. Then the solution $ (\hat{\vh},\hat{\vm}) $  of the LP in \eqref{eq:LP} lies in the set\footnote{ For a set $\setA \subset \R^m$, and a vector $\va \in \R^m$, we define by $\va \oplus \setA$, a set obtained by incrementing every element of $\setA$ by $\va$. } $(\tilde{\vh},\tilde{\vm})\oplus \setS$ with probability at least $1-\mathrm{e}^{-2Lt^2}$. 
\end{lemma}

Proof of this lemma is based on small ball method developed in \cite{koltchinskii2015bounding,mendelson2014learning} and further studied in  \cite{lecue2018regularization,lecue2017regularization}. The proof is mainly repeated using the argument in \cite{bahmani2017anchored}, and is provided in the supplementary material for completeness. We now state the main theorem for linear program \eqref{eq:LP}. The theorems states that if $L \geq 	C_t(S_1+S_2)\log^2(K+N)$, then the minimizer of the linear program \eqref{eq:LP} is in the set $(\tilde{\vh},\tilde{\vm}) \oplus \setS$ with high probability.

\begin{theorem}[Exact recovery]\label{thm:Main}
	Suppose we observe pointwise product of two vectors $\mB\ho$, and $\mC\mo$ through a bilinear measurement model in \eqref{eq:measurements}, where $\mB$, and $\mC$ are standard Gaussian random matrices. Then the linear program \eqref{eq:LP} recovers $(\hat{\vh},\hat{\vm}) \in (\tilde{\vh},\tilde{\vm}) \oplus \setS$ with probability at least $1-\mathrm{e}^{-(1/2)Lt^2}$ whenever $L \geq 	C_t(S_1+S_2)\log^2(K+N)$, where $C_t$ is a constant that depends on $t \geq 0$. 
\end{theorem} 

In light of Lemma \ref{lem:Mendelson}, the proof of Theorem \ref{thm:Main} reduces to computing the Rademacher complexity $\mathfrak{C}(\mathcal{D})$ defined in \eqref{eq:Rademacher-Complexity}, and the tail probability estimate $\mathfrak{p}_{\tau}(\mathcal{D})$ defined in \eqref{eq:Tail-Prob} of the set of descent directions $\mathcal{D}$ defined in \eqref{eq:setD}. These quantities are computed in the Supplementary material. The proof of Theorem \ref{thm:Noiseless_Main} follows by applying Lemma \ref{lem:feas} to Theorem \ref{thm:Main}.

\bibliographystyle{abbrv}

\newpage

\appendix
\section{Supplementary material}
{
\subsection{Proof of Lemma \ref{lem:feas}:} 
We first show that the feasible set of \eqref{eq:BH} is contained in the feasible set of \eqref{eq:LP}. We do this by using the fact that a convex set with a smooth boundary is contained in the halfspace defined by the tangent hyperplane at any point of the boundary of the convex set. Note that \eqref{eq:BH} is equivalent to the formulation
\begin{align*}
	\minimize_{\vh\in\R^K, \vm \in \R^N} \ \|\vh\|_1+\|\vm\|_1  \text{ subject to } &\yl \blt \vh \clt \vm  \geq \yl^2\\[-.5em]
	&t_\l \blt \vh \geq 0, ~ \l = 1, \ldots, L.
	\end{align*}
Consider a point $(\tilde{\wl} ,\tilde{\xl})$ on the boundary of the convex set defined by the constraints above and observe that 
	\begin{align}
	\left \{ (\wl, \xl) \in \R^2 \bigg | \begin{matrix}\yl \wl \xl \geq \yl^2 \\ \sign(\wl) = t_\l \end{matrix} \right \} \subseteq \left \{ (\wl, \xl) \in \R^2 \bigg | \begin{pmatrix}
	\yl \tilde{\xl} \\ \yl \tilde{\wl} \end{pmatrix} \cdot \begin{pmatrix}\wl - \tilde{\wl} \\ \xl - \tilde{\xl} \end{pmatrix} \geq 0 \right \}.
	\end{align}
	Plugging in $\wl = \blt \vh$ and $\xl = \clt \vm$, we have that any feasible $(\vh, \vm)$ satisfies 
	\[\yl \clt \tilde{\vm} \blt \vh + \yl \blt \tilde{\vh} \clt \vm \geq 2 \yl^2, \quad \ell = 1, \ldots, L,\]
	which implies $s_\l(\blt \vh \clt \tilde{\vm}  + \blt \tilde{\vh} \clt \vm) \geq 2 |\yl|$ for all $\ell$. So, the feasible set of \eqref{eq:LP} contains the feasible set of \eqref{eq:BH}. 
	
	Lastly, note that among all points $(\vh,\vm) \in (\tildeh,\tildem)\oplus S$, only $(\tildeh,\tildem)$ is feasible in \eqref{eq:BH}. So, if $(\tildeh,\tildem)$ solves \eqref{eq:LP} then $(\tildeh,\tildem)$ solves \eqref{eq:BH}.\qed
}
\subsection{Proof of Lemma \ref{lem:Mendelson}:} 
Define a one-sided loss function: 
\[
\mathcal{L}(\vh,\vm) := \tfrac{1}{L}\sum_{\ell = 1}^L \Big[2|\yl|-\sl\clt\tilde{\vm}\blt\vh-\sl\blt\tilde{\vh}\clt\vm\Big]_+,
\]
where $(\cdot)_+$ denotes the positive side. The LP in \eqref{eq:LP} can now be equivalently expressed as
\begin{align}\label{eq:LP-noisy}
&(\hat{\vh},\hat{\vm}):= \underset{(\vh,\vm) \in \R^{K+N}}{\text{argmin}} ~ \|\vh\|_1 + \|\vm\|_1  ~~\text{subject to}~~ \mathcal{L}(\vh,\vm)\leq 0. 
\end{align}
We want to show that there is no feasible descent direction $(\dh,\dm) \in \mathcal{D}$ around the true solution $(\tildeh,\tildem)$. Since  $(\dh,\dm)$ is a feasible perturbation from the proposed optimal $(\tildeh,\tildem)$, we have from \eqref{eq:LP-noisy}
\begin{align}\label{eq:interim}
\mathcal{L}(\tildeh+\dh,\tildem+\dm) \leq 0.
\end{align}
We begin by expanding the loss function $\mathcal{L}(\tildeh+\dh,\tildem+\dm)$ below
\begin{align}\label{eq:interim-eq1}
\mathcal{L}(\tildeh+\dh,\tildem+\dm) &= \tfrac{1}{L}\sum_{\ell=1}^L \big[\sl(2\yl-\blt\tildeh\clt(\tildem+\dm)-\clt\tildem\blt(\tildeh+\dh)\big]_+\notag\\
 &\geq\tfrac{1}{L}\sum_{\ell=1}^L\big[-\sl\blt\tildeh\clt\dm-\sl\clt\tildem\blt\dh\big]_+.
\end{align}
Let $\psi_t(s) := (s)_+-(s-t)_+$. Using the fact that $\psi_t(s) \leq (s)_+$, and that for every $\alpha, t \geq 0$, and $s \in \R$, $\psi_{\alpha t}(s) = t\psi_{\alpha}(\frac{s}{t})$, we have
\begin{align}\label{eq:interim-eq2}
&\tfrac{1}{L}\sum_{\ell=1}^L\big[-\sl\blt\tildeh\clt\dm-\sl\clt\tildem\blt\dh\big]_+ \geq \tfrac{1}{L}\sum_{\ell=1}^L\psi_{\tau \|(\dh,\dm)\|_2}\left(-\sl\blt\tildeh\clt\dm-\sl\clt\tildem\blt\dh\right) \notag\\
&\quad = \|(\dh,\dm)\|_2\cdot \tfrac{1}{L}\sum_{\ell=1}^L\psi_{\tau}\left(-\sl\left\<(\clt\tildem\bl,\blt\tildeh\cl),\tfrac{(\dh,\dm)}{\|(\dh,\dm)\|_2}\right\>\right)\notag\\
&\quad = \|(\dh,\dm)\|_2 \Bigg[\tfrac{1}{L}\sum_{\ell=1}^L\E\psi_{\tau}\left(-\sl\left\<(\clt\tildem\bl,\blt\tildeh\cl),\tfrac{(\dh,\dm)}{\|(\dh,\dm)\|_2}\right\>\right) - \notag\\
&  \tfrac{1}{L}\sum_{\ell=1}^L\bigg(\E\psi_{\tau}\left(-\sl\left\<(\clt\tildem\bl,\blt\tildeh\cl),\tfrac{(\dh,\dm)}{\|(\dh,\dm)\|_2}\right\>\right) - \psi_{\tau}\left(-\sl\left\<(\clt\tildem\bl,\blt\tildeh\cl),\tfrac{(\dh,\dm)}{\|(\dh,\dm)\|_2}\right\>\right)\bigg)\Bigg].
\end{align}
 The proof mainly relies on lower bounding the right hand side above uniformly over all $(\dh,\dm) \in \setD$. To this end, define a centered random process $\mathcal{R}(\mB,\mC)$ as follows
\begin{align*}
\mathcal{R}(\mB,\mC)&:=\sup_{(\dh,\dm)\in \setD}\tfrac{1}{L}\sum_{\ell=1}^L\bigg[\E\psi_{\tau}\bigg(-\sl\left\<(\clt\tildem\bl,\blt\tildeh\cl),\tfrac{(\dh,\dm)}{\|(\dh,\dm)\|_2}\right\>\bigg)\\
&\quad\qquad\qquad\qquad\qquad\qquad\qquad  -\psi_{\tau}\bigg(-\sl\left\<(\clt\tildem\bl,\blt\tildeh\cl),\tfrac{(\dh,\dm)}{\|(\dh,\dm)\|_2}\right\>\bigg)\bigg],
\end{align*}
and an application of bounded difference inequality \cite{mcdiarmid1989method} yields that $\mathcal{R}(\mB,\mC) \leq \E \mathcal{R}(\mB,\mC) + t/\sqrt{L}$  with probability at least $1-\mathrm{e}^{-2Lt^2/\tau^2}$. It remains to evaluate $\E  \mathcal{R}(\mB,\mC)$, which after using a simple symmetrization inequality \cite{van1997weak} yields 
\begin{align}
&\E \setR(\mB,\mC) \leq 2\E \sup_{(\dh,\dm)\in \setD\cap \setB}\tfrac{1}{L}\sum_{\ell=1}^L\varepsilon_\ell \psi_{\tau}\bigg(-\sl\left\<(\clt\tildem\bl,\blt\tildeh\cl),\tfrac{(\dh,\dm)}{\|(\dh,\dm)\|_2}\right\>\bigg),
\end{align}
where $\varepsilon_1, \varepsilon_2, \ldots, \varepsilon_L$ are independent Rademacher random variables. Using the fact that $\psi_t(s)$ is a contraction: $|\psi_t(\alpha_1)-\psi_t(\alpha_2)| \leq |\alpha_1-\alpha_2|$ for all $\alpha_1, \alpha_2 \in \R$, we have from the Rademacher contraction inequality \cite{ledoux2013probability} that 
\begin{align}\label{eq:random-process}
\E \sup_{(\dh,\dm)\in \setD}&\tfrac{1}{L}\sum_{\ell=1}^L\varepsilon_\ell \psi_{\tau}\bigg(-\sl\left\<(\clt\tildem\bl,\blt\tildeh\cl),\tfrac{(\dh,\dm)}{\|(\dh,\dm)\|_2}\right\>\bigg) \notag\\
& \qquad\qquad\leq \E \sup_{(\dh,\dm)\in \setD}\tfrac{1}{L}\sum_{\ell=1}^L-\varepsilon_\ell \sl\left\<(\clt\tildem\bl,\blt\tildeh\cl),\tfrac{(\dh,\dm)}{\|(\dh,\dm)\|_2}\right\>\notag\\
& \qquad\qquad = \E \sup_{(\dh,\dm)\in \setD}\tfrac{1}{L}\sum_{\ell=1}^L\varepsilon_\ell \left\<(\clt\tildem\bl,\blt\tildeh\cl),\tfrac{(\dh,\dm)}{\|(\dh,\dm)\|_2}\right\>,
\end{align}
where the last equality is the result of the fact that multiplying Rademacher random variables with signs does not change the distribution. In addition, using the facts that $t\mathbf{1}(s\geq t) \leq \psi_t(s)$, and that random vectors $\{(\clt\tildem\bl,\blt\tildeh\cl)\}_{\ell=1}^L$ are identically distributed and the distribution is symmetric, it follows 
\begin{align}\label{eq:tail-prob}
&\tau\PP\left(-\sl\left\<(\clt\tildem\bl,\blt\tildeh\cl),\tfrac{(\dh,\dm)}{\|(\dh,\dm)\|_2}\right\>\geq \tau\right) = \tau\PP\left(\left\<(\clt\tildem\bl,\blt\tildeh\cl),\tfrac{(\dh,\dm)}{\|(\dh,\dm)\|_2}\right\>\geq \tau\right) \notag\\
&= \tau\E \left[ \mathbf{1}\left(\left\<(\clt\tildem\bl,\blt\tildeh\cl),\tfrac{(\dh,\dm)}{\|(\dh,\dm)\|_2}\right\>\geq \tau\right)\right] \leq \E \psi_{\tau}\left( \left\<(\clt\tildem\bl,\blt\tildeh\cl),\tfrac{(\dh,\dm)}{\|(\dh,\dm)\|_2}\right\>\right). 
\end{align}
Plugging \eqref{eq:tail-prob}, and \eqref{eq:random-process} in \eqref{eq:interim-eq2}, we have 
\begin{align*}
&\tfrac{1}{L}\sum_{\ell=1}^L\Big[-\sl\left\<(\clt\tildem\bl,\blt\tildeh\cl),\tfrac{(\dh,\dm)}{\|(\dh,\dm)\|_2}\right\>\Big]_+ \geq \\
&\qquad\qquad\qquad \tau\|(\dh,\dm)\|_2\PP\left(\left\<(\clt\tildem\bl,\blt\tildeh\cl),\tfrac{(\dh,\dm)}{\|(\dh,\dm)\|_2}\right\> \geq \tau\right) \\
&\qquad\qquad -2\|(\dh,\dm)\|_2 \Big(\E \sup_{(\dh,\dm)\in \mathcal{D}}\tfrac{1}{L}\sum_{\ell=1}^L\varepsilon_\ell\left\<(\clt\tildem\bl,\blt\tildeh\cl),\tfrac{(\dh,\dm)}{\|(\dh,\dm)\|_2}\right\>-\tfrac{t}{\sqrt{L}}\Big)
\end{align*}
Combining this with \eqref{eq:interim} and \eqref{eq:interim-eq1}, we obtain the final result 
\begin{align*}
&\tau\|(\dh,\dm)\|_2\Bigg[\PP\left(\left\<(\clt\tildem\bl,\blt\tildeh\cl),\tfrac{(\dh,\dm)}{\|(\dh,\dm)\|_2}\right\> \geq \tau\right) \\
& \qquad \qquad  -2 \left(\E \sup_{(\dh,\dm)\in \mathcal{D}}\tfrac{1}{L}\sum_{\ell=1}^L\varepsilon_\ell\left\<(\clt\tildem\bl,\blt\tildeh\cl),\tfrac{(\dh,\dm)}{\|(\dh,\dm)\|_2}\right\>-\tfrac{t}{\sqrt{L}}\right)\Bigg]\leq 0.
\end{align*}
Using the definitions in \eqref{eq:complexity-calc}, and \eqref{eq:Tail-Prob}, we can write 
\begin{align*}
\|(\dh,\dm)\|_2 \left(\tau \mathfrak{p}_{\tau}(\mathcal{D})- \frac{(2\mathfrak{C}(\mathcal{D}) + t)}{\sqrt{L}}\right) \leq 0.
\end{align*}
It is clear that  choosing $L \geq \left( \frac{2\mathfrak{C}(\mathcal{D})+t}{\tau \mathfrak{p}_\tau(\mathcal{D})}\right)^2$ implies 
\begin{align*}
\|(\dh,\dm)\|_2 \leq 0,
\end{align*}
which directly means that $(\dh,\dm) = (0,0)$. Recall that $\setS \subset \setN$, and $\setD \perp \setN$, where $\setS$ is defined in \eqref{eq:setS}, this implies that the minimizer $(\hat{\vh},\hat{\vm})$ of the LP \eqref{eq:LP} resides in the set $(\tildeh,\tildem)\oplus \setS$. This completes the proof of Lemma \ref{lem:Mendelson}. 


\subsection{Proof of Theorem \ref{thm:Main}:}
	In light of Lemma \ref{lem:Mendelson}, the proof of Theorem \ref{thm:Main} comes down to computing the Rademacher complexity $\mathfrak{C}(\mathcal{D})$ defined in \eqref{eq:Rademacher-Complexity}, and the tail probability estimate $\mathfrak{p}_{\tau}(\mathcal{D})$ defined in \eqref{eq:Tail-Prob} of the set of descent directions $\mathcal{D}$ defined in \eqref{eq:setD}. 

	{\textbf{Upper Bound on Rademacher Complexity:}} We will start by evaluating $\mathfrak{C}(\mathcal{D})$ 
	\begin{align}\label{eq:complexity-calc}
	\mathfrak{C}(\mathcal{D}) &= \E \sup_{(\dh,\dm) \in \mathcal{D}} \tfrac{1}{\sqrt{L}} \sum_{\ell=1}^L \varepsilon_\ell \left\<(\clt\tildem\bl,\blt\tildeh\cl),\tfrac{(\dh,\dm)}{\|(\dh,\dm)\|_2}\right\>\notag\\
	&\leq \E \left\|\tfrac{1}{\sqrt{L}} \sum_{\ell=1}^L \varepsilon_\ell \left(\clt\tilde{\vm} \bl\vert_{\Gamma_h},\blt \tilde{\vh} \cl\vert_{\Gamma_m}\right)\right\|_2 \cdot \sup_{(\dh,\dm) \in \mathcal{D}}\left\|\tfrac{\left(\dh_{\Gamma_h},\dm_{\Gamma_m}\right)}{\|(\dh,\dm)\|_2}\right\|_2   
	\notag \\
	& \quad + \E  \left\|\tfrac{1}{\sqrt{L}} \sum_{\ell=1}^L \varepsilon_\ell \left(\clt\tilde{\vm} \bl\vert_{\Gamma^c_h},\blt \tilde{\vh} \cl\vert_{\Gamma^c_m}\right)\right\|_\infty\cdot  \sup_{(\dh,\dm) \in \mathcal{D}}\left\|\tfrac{\left(\dh_{\Gamma_h^c},\dm_{\Gamma_m^c}\right)}{\|(\dh,\dm)\|_2}\right\|_1.
	\end{align}
	First note that on set $\setD$ \eqref{eq:setD}, we have 
	\begin{align*}
	\left\|\tfrac{\big(\dh_{\Gamma_h^c},\dm_{\Gamma_m^c}\big)}{\|(\dh,\dm)\|_2}\right\|_1 \leq \sqrt{S_1+S_2} \left\|\tfrac{\left(\dh_{\Gamma_h},\dm_{\Gamma_m}\right)}{\|(\dh,\dm)\|_2}\right\|_2  \leq \sqrt{S_1+S_2}. 
	\end{align*}
   As for the remaining terms, we begin by writing
	\begin{align*}
	&\E 	\left\|\tfrac{1}{\sqrt{L}} \sum_{\ell=1}^L \varepsilon_{\ell} \left(\clt\tilde{\vm} \bl\vert_{\Gamma_h},\blt \tilde{\vh} \cl\vert_{\Gamma_m}\right)\right\|_2 \leq \sqrt{\E 	\left\|\tfrac{1}{\sqrt{L}} \sum_{\ell=1}^L \varepsilon_{\ell} \left(\clt\tilde{\vm} \bl\vert_{\Gamma_h},\blt \tilde{\vh} \cl\vert_{\Gamma_m}\right)\right\|_2^2}\\
	& \qquad\qquad \qquad = \sqrt{\tfrac{1}{L}\sum_{\ell=1}^L \E \left(|\clt\tilde{\vm}|^2 \|\blt\vert_{\Gamma_h}\|_2^2 +|\bl\tilde{\vh}|^2 \|\cl\vert_{\Gamma_m}\|_2^2\right) }\\ 
	&\qquad\qquad\qquad= \sqrt{\|\tilde{\vm}\|_2^2 S_1+\|\tilde{\vh}\|_2^2 S_2},
   \end{align*}
   and the second term in \eqref{eq:complexity-calc} is  
   \begin{align*}
   &\E \left\|\tfrac{1}{\sqrt{L}} \sum_{\ell=1}^L \varepsilon_\ell(\blt\tilde{\vh}\cl\vert_{\Gamma_m^c},\clt\tilde{\vm}\bl\vert_{\Gamma_h^c})\right\|_\infty \leq \sqrt{  \E \left\|\tfrac{1}{\sqrt{L}} \sum_{\ell=1}^L \varepsilon_\ell(\blt\tilde{\vh}\cl\vert_{\Gamma_m^c},\clt\tilde{\vm}\bl\vert_{\Gamma_h^c})\right\|^2_\infty }\\
   & \qquad\qquad\qquad \leq \sqrt{2e \log (K+N)\cdot \tfrac{1}{L}\sum_{\ell=1}^L\E \max\left\{|\clt\tilde{\vm}|^2 \|\bl\vert_{\Gamma_h^c}\|_{\infty}^2,|\blt\tilde{\vh}|^2 \|\cl\vert_{\Gamma_m^c}\|_\infty^2\right\}}\\
   & \qquad\qquad\qquad\leq \sqrt{2e\log (K+N)\E\max\{|\vb^{\intercal}\tilde{\vh}|^2\|\vc\vert_{\Gamma_m^c}\|_{\infty}^2, |\vc^{\intercal}\tilde{\vm}|^2\|\vb\vert_{\Gamma_h^c}\|_{\infty}^2\}}\\
  &\qquad\qquad\qquad\leq C\sqrt{\max\{\|\tilde{\vh}\|_2^2,\|\tilde{\vm}\|_2^2\}\log^2(K+N)}, 
   \end{align*}
   where the second inequality by the application of Lemma 5.2.2 in \cite{akritas2016topics}, 
and the final equality is due to the fact that $\|\vc\vert_{\Gamma_m^c}\|_{\infty}^2$, and $\|\vb\vert_{\Gamma_h^c}\|_{\infty}^2$  are subexponential and using Lemma 3 in \cite{van2013bernstein}.

Plugging the bounds above back in \eqref{eq:complexity-calc}, we obtain the upper bound on the Rademacher complexity given below
   \begin{align}\label{eq:complexity-upperbound}
   \mathfrak{C}(\mathcal{D}) &\leq C\sqrt{\big(\|\tilde{\vm}\|_2^2+\|\tilde{\vh}\|_2^2\big)  (S_1+ S_2)\log^2(K+N)}.
   \end{align}
{\textbf{Tail Probability: }}To apply the result in Lemma \ref{lem:Mendelson}, we also need to evaluate 
\begin{align}\label{eq:tail-probability}
\mathfrak{p}_{\tau}(\mathcal{D})= \inf_{(\dh,\dm) \in \setD}
   \PP\left(\left\<(\clt\tildem\bl,\blt\tildeh\cl),\tfrac{(\dh,\dm)}{\|(\dh,\dm)\|_2}\right\>\geq \tau\right). 
  \end{align}
   It suffice to estimate the probability $\PP(|\blt\tilde{\vh}\clt\dm+\blt\dh\clt\tilde{\vm}| \geq \tau)$, which using Paley-Zygmund inequality implies 
   \begin{align*}
   & \PP\left(|\blt\tilde{\vh}\clt\dm+\blt\dh\clt\tilde{\vm}|^2 \geq \frac{1}{2}\E |\blt\tildeh\clt\dm+\blt\dh\clt\tildem|^2 \right) \\
   & \qquad\qquad \qquad \geq \frac{1}{4}\cdot \frac{(\E |\blt\tildeh\clt\dm+\blt\dh\clt\tildem|^2 )^2}{\E |\blt\tilde{\vh}\clt\dm+\blt\dh\clt\tilde{\vm}|^4}.
   \end{align*}
   Using norm equivalence of Gaussian random variables, we know that $(\E |\blt\tildeh\clt\dm+\blt\dh\clt\tildem|^4)^{1/4} \leq c (\E |\blt\tildeh\clt\dm+\blt\dh\clt\tildem|^2)^{1/2}$, this implies that 
   \begin{align}
   \PP(|\blt\tildeh\clt\dm+\blt\dh\clt\tildem|^2 \geq \frac{1}{2}\E |\blt\tildeh\clt\dm+\blt\dh\clt\tildem|^2 ) \geq \frac{1}{4}\cdot \frac{1}{c^4}.
   \end{align}
 Finally, a simple calculation shows that $\E |\blt\tildeh\clt\dm+\blt\dh\clt\tildem|^2 \geq c(\|\tildeh\|_2^2\|\dm\|_2^2 +\|\tildem\|_2^2\|\dh\|_2^2)$ for an absolute constant $c$. 
\begin{align*}
\E |\blt\tildeh\clt\dm+\blt\dh\clt\tildem|^2 &= \E_{\vb}\E_{\vc} \tildeh^\top \bl \blt\tildeh\dm^\top \cl \clt\dm + \dh^\top \bl\blt\dh\tildem^\top\cl \clt\tildem \\ &
\qquad \quad + 2\E_{\vb}\E_{\vc} \dh^\top\bl\blt\tildeh\dm^\top \cl\clt\tildem \\ & = \E_{\vb} \|\dm\|^2 \tildeh^\top \bl \blt\tildeh + \|\tildem\|^2 \dh^\top \bl\blt\dh + 2\dm^\top \tildem \dh^\top\bl\blt\tildeh\\ 
&=\|\dm\|^2\|\tildeh\|^2 + \|\tildem\|^2\|\dh\|^2 + 2 \dm^\top \tildem \dh^\top \tildeh\\
& =\|\dm\|^2\|\tildeh\|^2 + \|\tildem\|^2\|\dh\|^2 + 2 (\dh^\top \tildeh)^2\\
& \geq \|\dm\|^2\|\tildeh\|^2 + \|\tildem\|^2\|\dh\|^2,
\end{align*}
where the last equality follows using the fact $(\dh,\dm) \in \setD \subset \setN_\perp$, and hence $\setD \perp \setN$, which implies that $\dh^\top\tildeh = \dm^\top\tildem.$ Normalizing by $\|(\dh,\dm)\|_2$, and comparing with \eqref{eq:tail-probability} directly shows that $\tau^2 = (\|\tildeh\|_2^2\|\dm\|_2^2 +\|\tildem\|_2^2\|\dh\|_2^2)$, and $\mathfrak{p}_{\tau}(\mathcal{D}) = 0.5/c^4$. Plugging these results and the Rademacher complexing bound in \eqref{eq:complexity-upperbound}, in Lemma \ref{lem:Mendelson} proves Theorem \ref{thm:Main}. 	
\qed

\subsection{Evaluation of the Projection Operator}\label{sec:proj}

Given a point $(\vx' ,\vw',\bxi')\in\mathbb{R}^{3L}$, in this section we focus on deriving a closed-form expression for $\mbox{proj}_{\mathcal{C}}\left( (\vx' ,\vw',\bxi')\right)$, where 
\[C = \left\{(\vx,\vw,\vxi)\in\R^{3L}|\ s_\ell(\xi_\ell+x_\ell)w_\ell\geq |y_\ell|,\ t_\ell w_\ell\geq 0,\ \l = 1,\dots, L\right\}\]
 is the convex feasible set of $\eqref{eq:Rl1BH}$. It is straightforward to see that the resulting projection program decouples into $L$ convex programs in $\mathbb{R}^3$ as
\begin{equation}\label{projeq}\argmin_{x\in\R,w\in\R,\xi\in\R}~~ \frac{1}{2}\left\|\begin{pmatrix}x\\ w\\ \xi \end{pmatrix} -  \begin{pmatrix} x_\ell'\\ w_\ell'\\ \xi_\ell'  \end{pmatrix}\right\|_2^2 ~~s.t.~~ |y_\ell| - s_\ell xw - s_\ell \xi  w  \leq 0,\quad  -t_\ell w \leq 0.
\end{equation}
Throughout this derivation we assume that $|y_\ell|>0$ (derivation of the projection for the case $y_\ell$ is easy) and as a result of which the second constraint $-t_\ell w \leq 0$ is never active (because then $w=0$ and the first constraint requires that $|y_\ell|\leq 0$). We also consistently use the fact that $t_\ell$ and $s_\ell$ are signs and nonzero. 

Forming the Lagrangian as
\[\mathcal{L}(x,w,\xi,\mu_1,\mu_2) = \frac{1}{2}\left\|\begin{pmatrix}x\\ w\\ \xi \end{pmatrix} -  \begin{pmatrix} x_\ell'\\ w_\ell'\\ \xi_\ell'  \end{pmatrix}\right\|_2^2 + \mu_1\left( |y_\ell| - s_\ell xw - s_\ell \xi w \right ) - \mu_2\left( t_\ell w \right),
\]
along with the primal constraints, the KKT optimality conditions are 
\begin{align}\label{e5}
\frac{\partial \mathcal{L}}{\partial x} = x-x_\ell' - \mu_1s_\ell w   &=0,\\ \label{e6}
\frac{\partial \mathcal{L}}{\partial w} = w - w_\ell' - \mu_1s_\ell  x - \mu_1s_\ell  \xi -\mu_2 t_\ell  &=0,\\ 
 \label{e7}
\frac{\partial \mathcal{L}}{\partial \xi} = \xi-\xi_\ell' - \mu_1 s_\ell w&=0,\\
\label{e8}
\mu_1\geq 0, \quad \mu_1\left( |y_\ell| - s_\ell xw - s_\ell \xi w \right ) &=0,\\
 \label{e9} \mu_2 \geq 0, \quad \mu_2\left( t_\ell w \right)&=0.
\end{align}
We now proceed with the possible cases.

\textbf{Case 1.} $\mu_1=\mu_2=0$:\\
In this case we have $(x,w,\xi)=(x_\ell',w_\ell',\xi_\ell')$ and this result would only be acceptable when $|y_\ell| - s_\ell x_\ell'w_\ell' - s_\ell \xi_\ell' w_\ell' \leq 0$ and $t_\ell w_\ell'\geq 0$.

\textbf{Case 2.} $\mu_1=0$, $t_\ell w =0$:\\ 
In this case the first feasibility constraint of \eqref{projeq} requires that $|y_\ell|\leq 0$, which is not possible when $|y_\ell|>0$.

\textbf{Case 3.} $|y_\ell| - s_\ell xw - s_\ell \xi w = 0$, $t_\ell w =0$:\\ 
Similar to the previous case, this cannot happen when $|y_\ell|>0$.

\textbf{Case 4.} $\mu_2=0$, $|y_\ell| - s_\ell xw - s_\ell \xi w  =0$:\\ 
In this case we have
\[|y_\ell| = s_\ell xw + s_\ell \xi w.
\] 
Now combining this observation with \eqref{e5} and \eqref{e7} yields 
\begin{align}\label{e18}
|y_\ell| = s_\ell \left(   x_\ell' + \mu_1s_\ell w \right) w + s_\ell \left(   \xi_\ell' + \mu_1s_\ell w \right) w,
\end{align}
and therefore
\begin{align}\label{e17}
\mu_1 = \frac{|y_\ell| - s_\ell \left(   x_\ell' + \xi_\ell' \right)w}{2w^2}.
\end{align}
Similarly, \eqref{e6} yields 
\begin{equation}\label{e19}
w=w_\ell'+ \mu_1 s_\ell \left(   x_\ell' + \mu_1s_\ell w \right)  + \mu_1 s_\ell \left(   \xi_\ell' + \mu_1s_\ell w \right).
\end{equation}
Knowing that $w\neq 0$, $\mu_1$ can be eliminated between \eqref{e18} and \eqref{e19} to generate the following forth order polynomial equation in terms of $w$:
\begin{align*}
2w^4-2w_\ell'w^3 +s_\ell|y_\ell|\left(x_\ell'+\xi_\ell' \right)w- y_\ell^2=0.
\end{align*}
After solving this 4-th order polynomial equation (e.g., the root command in MATLAB) we pick the real root $w$ which obeys
\begin{align}\label{eqconsts}
t_\ell  w\geq 0, \qquad |y_\ell| - s_\ell \left(   x_\ell' + \xi_\ell' \right) w\geq 0.
\end{align}
Note that the second inequality in \eqref{eqconsts} warrants nonnegative values for $\mu_1$ thanks to \eqref{e17}. After picking the right root, we can explicitly obtain $\mu_1$ using \eqref{e19} and calculate the solutions $x$ and $\xi$ using \eqref{e5} and \eqref{e7}.
Technically, in using the ADMM scheme for each $\ell$ we solve a forth-order polynomial equation and find the projection.
\end{document}